# A nested hierarchy of second order upper bounds on system failure probability

Sourangshu Ghosh

Graduate Student

Department of Civil Engineering

Indian Institute of Technology Kharagpur

Kharagpur 721303, WB, India

Email: sourangshu@iitkgp.ac.in

Baidurya Bhattacharya

Professor

Department of Civil Engineering

Indian Institute of Technology Kharagpur

Kharagpur 721303, WB, India

Corresponding author. Email: baidurya@civil.iitkgp.ac.in

**ABSTRACT**

For a coherent, binary system made up of binary elements, the exact failure probability requires knowledge of statistical dependence of all orders among the minimal cut sets. Since dependence among the cut sets beyond the second order is generally difficult to obtain, second order bounds on system failure probability have practical value. The upper bound is conservative by definition and can be adopted in reliability based decision making. In this paper we propose a new hierarchy of $m$-level second order upper bounds, $B_m$ : the well-known Kounias-Vanmarcke-Hunter-Ditlevsen (KVHD) bound - the current standard for upper bounds using second order joint probabilities - turns out to be the weakest member of this family ($m = 1$). We prove that $B_m$ is non-increasing with level $m$ in every ordering of the cut sets, and derive conditions under which $B_{m+1}$ is strictly less than $B_m$ for any $m$ and any ordering. We also derive conditions under which the optimal level $m$ bound is strictly less than the optimal level $m + 1$ bound, and show that this improvement asymptotically achieves a probability of 1 as long as the second order joint probabilities are only constrained by the pair of corresponding first order probabilities. Numerical examples show that our second order upper bounds can yield tighter values than previously achieved and in every case exhibit considerable less scatter across the entire $n$! orderings of the cut sets compared to KVHD bounds. Our results therefore may lead to more efficient identification of the optimal upper bound when coupled with existing linear programming and tree search based approaches.

## 1. Introduction

For a binary system made up of binary elements, the system failure event can be described as the union of its minimal cut sets:

$$F_{sys} = \bigcup_{i=1}^{n} C_i \tag{1}$$

Each minimal cut set, $C_i$, is a parallel arrangement of its constituent elements:

$$C_i = F_{i_1} \cap F_{i_2} \cap ... \cap F_{i_{max}}, \quad i = 1,...,n \tag{2}$$

where $F_i = \{X_i = 0\}$, $i = 1,...,n_{el}$ is the failure of the $i^{th}$ binary element with

$$X_i = \begin{cases} 0 \text{ if element } i \text{ is down} \\ 1 \text{ if element } i \text{ is up} \end{cases}, \quad i = 1,...,n_{el} \tag{3}$$

The minimal cut sets are generally not independent (nor are they disjoint) owing to (i) the presence of the same element failure event $F_j$ in more than one $C_i$'s, and (ii) possible mutual dependence among the $F_j$'s themselves. Hence, a central problem in reliability analysis is to estimate the union probability in Eq (1):

$$P[F_{sys}] = P\bigcup_{i=1}^{n} C_i = \sum_{\text{all } i} P_i - \sum_{\text{all } i,j; j<i} P_{ij} + \sum_{\text{all } i,j,k; k<j<i} P_{ijk} - ... \tag{4}$$

where $P_i = P[C_i], P_{ij} = P[C_i C_j], P_{ijk} = P[C_i C_j C_k]$, etc. In general, the evaluation of $P_i, P_{ij}, P_{ijk},...$ requires the joint probability information of the constituent element failure events $F_{i_j}$. If each cut set in (1) can be described by a limit state function $g_i$ such that $C_i = \{g_i < 0\}$ and $g_i$ is a linear combination of one or more jointly normal random variables, then an exact (numerical) evaluation of Eq (4) is possible with only the first order $P_i$'s and the second order $P_{ij}$'s; in every other case, higher order joint probabilities are required for evaluating the union probability.

Bonferroni [1] first introduced upper and lower bounds which are simple algebraic sums with alternating signs of the joint probabilities. As a matter of practical consideration however, joint probabilities beyond the second order are difficult to obtain, and hence bounds on the union probability based only on second order joint information have a practical appeal. Further, second order upper bounds, which are the subject of this paper, are again of practical interest in reliability analysis as they provide conservative estimates of system failure probability with limited data.

Upper and lower Bonferroni bounds have been the subject of considerable research since the 1950s. The first such (lower) bound was discovered by Chung and Erdos [2] which was also found independently by Whittle [3]. The first approximations to the union probability in structural reliability involved only first order probabilities: Freudenthal et al. [4] approximated it as the sum of first order failure probabilities $\sum P_i$. Cornell [5] proposed the lower bound to the union probability as $\max P_i$ and showed that for a coherent system the upper bound to the union probability is $1 - \prod_{i=1}^{n}[1 - P(F_i)]$. Using the Bonferroni inequalities, Kounias [6] obtained upper and lower bounds involving both first and second order probabilities; similar second order bounds were subsequently proposed by Vanmarcke [7], Hunter [8] and Ditlevsen [9]: We refer to these as KVHD bounds in this paper. For a structural system with normally distributed performance functions, Ahmed and Koo [10] showed that the upper and lower bounds of the resultant joint normal probability are narrower than KVHD bounds. Improvements using third or higher order joint probabilities to KVHD second order bounds were later proposed by Hohenbichler and Rackwitz [11], Ramachandran [12], Feng [13], Greig [14], Zhang [15], and Ramachandran [16]. Reliability bounds based on interval probability theory have been developed by Cui and Blockley [17], Qiu et al. [18], and Wang et al. [19, 20]. Recently, a new method using interval Monte Carlo method along with Linear Programming has been developed by Zhang et al. [21].

The second and higher order bounds discussed above depend on the ordering of the failure events and one would in principle need to compute the bounds for all $n!$ permutations of the minimal cut sets in order to obtain the sharpest bounds. This can be computationally expensive for large problems and researchers have looked for methods that do not require computing bounds for all orderings. Hailerpin [22] was the first to formulate the Boolean probability bounding problem as a linear programming (LP) problem and showed that Boole's method is similar to Fourier's elimination. Using the LP proposed by Hailerpin, Kounias and Marin [23] proposed second order upper and lower bounds using indicator random variables and LP. They showed that previously known bounds [2, 3, 6, 24-27] are particular cases of their bounds. They have also shown that that if the events are assumed exchangeable then their bounds are the best in a given class of bounds. Around the same time, Kwerel [28, 29] described the dual feasible bases of LP to obtain upper bounds on union probabilities based on first two binomial moments. Galambos [30] also found the same upper bound based on first two binomial moments using a different technique. A few years later, Galambos and Mucci [31] and Platz [32] developed bounds using LP that use higher binomial moments. Prekopa in his series of papers [33-36] formulated the Bonferri Inequalities of Dawson and Sankoff [24] as a linear programming problem, replaced the first and second order probabilities with the first $m$ binomial moments of the random variable and obtained sharper bounds.

Tree structures have also been used to search optimal bounds. Bukszár and Prékopa [37] introduced the idea of Cherry Trees which are special cases of chordal graph structure, and derived third order upper bounds to the union probability. Tomescu [38] generalized the Hunter Bound [8] and also proposed new lower bounds using the concept of hypertrees in the framework of uniform hypergraphs. Bukszár and Szántai [39] improved Tomescu's lower and upper bounds [38] by introducing the idea of hypercherry tree in the same ways as Bukszár and Prékopa [37] generalized the Hunter-Worsley [8, 40] bound. Boros and Veneziani [41] generalized the cherry tree bounds by using chordal graph structure which are graphs where every cycle of 4 or more vertices have a chord that connects two non-consecutive vertices of the cycle. This graph structure was further generalized by Dohmen [42, 43] to find a new set of lower bounds using chordal-sieve bounds.

For structural systems, Song and Der Kiureghian [44] showed for the first time that LP can be used to compute bounds given any available information on the component probabilities and that the LP based bounds were independent of the ordering of the components and produced the narrowest possible bounds. Subsequently, Der Kiureghian and Song [45] extended the formulation to complex systems having large number of cut and link sets and proposed multi-scale modeling of the decomposed system. Chang and Mori [46] developed a relaxed linear programming (RLP) bounds method while Chang et al. [47] derived bounds on failure probability of $k$-out-of-$n$ systems with the help of universal generating function and LP. Byun and Song [48] applied binary integer programming to tackle the problem of exponential rise in the number of variables in LP with system size. A recent overview of all these structural reliability estimation methods is available in Song, Kang, Lee and Chun [49].

A considerable amount of work over the past decades has focused exclusively on the lower bound. Although outside the scope of this paper, we summarize them for the sake of completeness. Prekopa and Gao [50] generalized the lower bounds developed by De Caen [51] and Kuai et al. [52] using additional information (third order joint probabilities). The Kuai et al. [52] lower bound was further improved by Yang et al. [53, 54]. A similar lower bound using only first and second order probabilities was also proposed much earlier by Gallot [25]. This bound was recently revisited by Feng et al. [55, 56] and Mao et al. [57]. They also showed that the Gallot Bound [25] is not necessarily sharper than the Kuai et al. [52] lower bound. The De Caen bound [51] was further improved by Cohen and Merhav [58]. Szántai [59] used variance reduction technique to improve previously discovered lower bounds.

The union probability bounding problem is a special case of probabilistic satisfiability problem [60]. The linear programming models are generally computationally very intensive and not polynomially computable [61]. Zemel [62], Jaumard et al. [63] and Georgakopoulos et al. [60] proposed column generation techniques to solve this problem. Nevertheless, column generation and quadratic binary optimization are similar algorithms and thus column generation method is an NP-hard optimization problem [63]. Deza and Laurent [64] showed that column generation is algorithmically is similar to the separation problem for the cut polytope. They developed upper and lower bounds by using inequalities for this correlation polytope. Boros and Hammer [65] further generalize these cut polytope bounds. This complexity, feasibility of the cut polytope problem has been also discussed by Kavvadias and Papadimitriou [66] and Veneziani [67].

In this article, we propose a new hierarchy of $m$ level of second order upper bounds, $B_m$, to the $n$-dimensional ($m < n$) union probability $P[F_{sys}]$: The well-known Kounias-Vanmarcke-Hunter-Ditlevsen (KVHD) second order upper bound [6-9] turns out to be the weakest member of this family ($m = 1$). The hierarchy of bounds is non-increasing with level $m$ in every ordering of the cut sets, and we derive conditions under which $B_{m+1}$ is strictly less than $B_m$ for any $m$ and any ordering. We also derive conditions under which the optimal level $m + 1$ bound is strictly less than the optimal level $m$ bound, and show that this improvement asymptotically achieves a probability of 1 as long as the second order joint probabilities are only constrained by the pair of corresponding first order probabilities. Numerical examples show that our second order upper bounds can yield tighter values than previously achieved and in every case our bounds exhibit considerable less scatter across the entire $n!$ orderings of the cut sets compared to KVHD bounds which are the current standard for upper bounds using second order joint probabilities. Our results therefore may lead to more efficient identification of the optimal upper bound when coupled with existing linear programming and tree search based approaches.

Before presenting the general form, we start with deriving the level 2 bound, and show that even for $m = 2$, our second order bound is less sensitive to the ordering of the cut sets, that it is at least as good as the KVHD bound in every case, and, under a very mild condition, is better than the KVHD upper bound in a given ordering. The level 2 upper bound is given in Eq (11) and the general level $m$ upper bound is given in Eq (26) below.

## 2. The level-2 second order bound

We list out the contribution of each additional cut set in the union by rewriting Eq (4) as:

$$
\begin{aligned}
P[F_{sys}] = &\; P_1 \\
& + P_2 - P_{12} \\
& + P_3 - P_{13} - P_{23} + P_{123} \\
& + P_4 - P_{14} - P_{24} - P_{34} + P_{124} + P_{134} + P_{234} - P_{1234} \\
& + P_5 - P_{15} - P_{25} - P_{35} - P_{45} + P_{125} + P_{135} + P_{145} + P_{235} + P_{245} + P_{345} - P_{1235} - P_{1245} - P_{1345} - P_{2345} + P_{12345} \\
& + P_6 - ...
\end{aligned}
\quad (5)
$$

From the third line onward, we can rewrite (5) as:

$$P[F_{sys}] = P_1$$
$$+ P_2 - P_{12}$$
$$+ P_3 - P(C_1C_3 \cup C_2C_3)$$
$$+ P_4 - P(C_1C_4 \cup C_2C_4 \cup C_3C_4) \quad (6)$$
$$+ P_5 - P(C_1C_5 \cup C_2C_5 \cup C_3C_5 \cup C_4C_5)$$
$$+ P_6 - ...$$

Since $P(A_1 \cup A_2 \cup ...) \geq \max P(A_i)$ for any collection of sets $A_1, A_2, ...$, we have:

$$P[F_{sys}] \leq P_1$$
$$+ P_2 - P_{12}$$
$$+ P_3 - \max(P_{13}, P_{23})$$
$$+ P_4 - \max(P_{14}, P_{24}, P_{34})$$
$$+ P_5 - \max(P_{15}, P_{25}, P_{35}, P_{45}) \quad (7)$$
$$+ P_6 - ...$$
$$= P_1 + P_2 - P_{12} + \sum_{i=3}^{n} \left[ P_i - \max_{1 \leq j < i} \{P_{ji}\} \right] = B_1$$

which is the well-known second order KVHD upper bound [6-9] mentioned above. In this paper we show that KVHD upper bound happens to be the first member of a family of hierarchical level-$m$ second order upper bounds, $B_m$, whose general form will be presented in Section 4. Before presenting the general form, we present the level 2 bound next.

We can obtain a better bound by going back to the third line onward in (6). Since $P(A_1 \cup A_2 \cup A_3 ...) \geq P(A_i \cup A_j), i, j = 1, 2, 3, ..., i \neq j$ for any collection of three or more sets $A_1, A_2, A_3, ...$, we have:

$$P[F_{sys}] \leq P_1$$
$$+ P_2 - P_{12}$$
$$+ P_3 - P(C_1C_3 \cup C_2C_3)$$
$$+ P_4 - \max \left[ P(C_1C_4 \cup C_2C_4), P(C_1C_4 \cup C_3C_4), P(C_2C_4 \cup C_3C_4) \right]$$
$$+ P_5 - \max \begin{bmatrix} P(C_1C_5 \cup C_2C_5), P(C_1C_5 \cup C_3C_5), P(C_1C_5 \cup C_4C_5), P(C_2C_5 \cup C_3C_5), \\ P(C_2C_5 \cup C_4C_5), P(C_3C_5 \cup C_4C_5) \end{bmatrix} \quad (8)$$
$$+ P_6 - ...$$
$$= P_1 + P_2 - P_{12} + \sum_{i=3}^{n} \left[ P_i - \max_{1 \leq j < l < i} \left[ P(C_jC_i \cup C_lC_i) \right] \right]$$

Let us look at any one argument within the max [ ] brackets in (8). The general form is:

$$P(C_j C_i \cup C_l C_i) = P_{ji} + P_{li} - P_{jli} \qquad (9)$$

Since $P_{jli} \leq P_{ji}, P_{jli} \leq P_{li}, P_{jli} \leq P_{lj}$ in all cases, we can write:

$$P(C_j C_i \cup C_l C_i) \geq P_{ji} + P_{li} - \min(P_{ji}, P_{li}, P_{lj}) \qquad (10)$$

which gives us a new upper bound:

$$P[F_{sys}] \leq P_1 + P_2 - P_{12} + \sum_{i=3}^{n} \left[ P_i - \max_{1 \leq j < l < i} \left\{ P_{ji} + P_{li} - \min(P_{ji}, P_{li}, P_{lj}) \right\} \right]^+ = B_2 \qquad (11)$$

We first show that this level 2 bound is at least as good as KVHD bound in every permutation of the index set, and then derive the condition under which $B_2$ is better than $B_1$ in a given permutation. Subsequently, we discuss under what conditions the best $B_2$ is better than the best $B_1$ over all permutation of the index set. We will also generalize the results as the number of cut sets (*n*) becomes large.

## 3. An improvement over KVHD bound

The proposed level 2 upper bound (11) is always less than or equal to the upper KVHD bound regardless of the ordering of events; further, if a rather mild condition is satisfied (which we term Condition 1 below), there are at least $2(n-3)!$ orderings where our bound is strictly less than KVHD. To show these we need the following results.

*Theorem* 1. In any ordering $(\pi)$ of the index set describing second order probabilities, the level 2 bound is less than or equal to the corresponding level 1 bound: $B_2(\pi) \leq B_1(\pi)$.

*Proof.* We prove the theorem by showing that for all quantities $P_{ji}, P_{li}$ and $P_{lj}$ such that $1 \leq j < l < i, 3 \leq i$ in ordering $(\pi)$, we must have

$$\max_{1 \leq j < l < i} \left\{ P_{ji} + P_{li} - \min(P_{ji}, P_{li}, P_{lj}) \right\} \geq \max_{1 \leq j < i} \left\{ P_{ji} \right\} \qquad (12)$$

For any three quantities *a*, *b* and *c* we can write:

$$b \geq \min(a,b,c)$$
$$b - \min(a,b,c) \geq 0 \qquad (13)$$

Adding *a* on both sides, we obtain

$$a + b - \min(a,b,c) \geq a \tag{14}$$

Without loss of generality, let us assign $a = P_{ji}, b = P_{li}, c = P_{lj}$. Taking the maximum on both sides of (14) over $1 \leq j < l < i$, $3 \leq i$ we arrive at (12). We now sum both sides of (12) from $i = 3$ to $n$ and subtract both sides from $P_1 + P_2 - P_{12} + \sum_{i=3}^{n} P_i$ to obtain:

$$P_1 + P_2 - P_{12} + \sum_{i=3}^{n} \left[ P_i - \max_{1 \leq j < l < i} \{ P_{ji} + P_{li} - \min(P_{ji}, P_{li}, P_{lj}) \} \right]$$
$$\leq P_1 + P_2 - P_{12} + \sum_{i=3}^{n} \left[ P_i - \max_{1 \leq j < i} \{ P_{ji} \} \right] \tag{15}$$

i.e., $B_2(\pi) \leq B_1(\pi)$

Hence, proved.

Since this holds for any ordering $(\pi)$ of the minimal cut sets $\{C_i\}$, i.e., for every permutation of the index set $\{1, 2, ..., n\}$, our bound (11) is at least as good as KVHD bound in Eq (7) for any given permutation of the cut sets. We now show that our bound is strictly better than KVHD under a rather mild condition, introduced next.

*Condition* 1: Given second order probabilities $P_{ij} = P_{ji}$, $i = 1,..., n-1, i < j \leq n$, in some ordering of the index set, there is one triplet $a,b,c$ (all distinct with $a,b < c$) for which the largest off-diagonal element above the diagonal in column $c$, $P_{ac} = \max_{i<c}(P_{ic})$, satisfies

$$P_{ac} = \max_{i<c}(P_{ic}) \geq P_{bc} > P_{ab} \tag{16}$$

*Theorem* 2: If a particular ordering of the index set of second order probabilities satisfies Condition 1, the level 2 bound is less than the level 1 bound in that ordering.

*Proof*: Since $P_{ac} \geq P_{bc} > P_{ab}$, we can write

$$P_{ac} + P_{bc} - \min(P_{ac}, P_{bc}, P_{ab}) > P_{ac} = \max_{i<c}(P_{ic}) \tag{17}$$

We have already proved (Theorem 1) that for any $1 \leq j < l < i$, $3 \leq i$

$$\max_{1 \leq j < l < i} \{ P_{ji} + P_{li} - \min(P_{ji}, P_{li}, P_{lj}) \} \geq \max_{1 \leq j < i} \{ P_{ji} \} \tag{18}$$

Summing both sides from $i = 3,..., n$, but $i \neq c$, we have

$$\sum_{\substack{\text{all } i \\ i \neq c}} \max_{1 \leq j < l < i} \{P_{ji} + P_{li} - \min(P_{ji}, P_{li}, P_{lj})\} \geq \sum_{\substack{\text{all } i \\ i \neq c}} \max_{1 \leq j < i} \{P_{ji}\} \tag{19}$$

Combining (17) with (19) and subtracting both sides from $P_1 + P_2 - P_{12} + \sum_{i=3}^{n} P_i$ we obtain:

$$P_1 + P_2 - P_{12} + \sum_{i=3}^{n} \left[ P_i - \max_{1 \leq j < l < i} \{P_{ji} + P_{li} - \min(P_{ji}, P_{li}, P_{lj})\} \right]$$
$$< P_1 + P_2 - P_{12} + \sum_{i=3}^{n} \left[ P_i - \max_{1 \leq j < i} \{P_{ji}\} \right] \tag{20}$$

i.e., $B_2 < B_1$ under Condition 1

Hence, proved.

If Condition 1 is satisfied for a certain $c$ in a given ordering of the index set $\{1, 2, ..., n\}$, it will be satisfied for a subset of other orderings of the index set as well. The minimum number of such orderings, $\sum_{j=0}^{c-3} \binom{c-3}{j}(j+2)!(n-3-j)!$, depends on the value of $c$ in (16), $3 \leq c \leq n$ where $j$ signifies the number of free columns (other than $a$ and $b$) to the left of the $c^{\text{th}}$ column: for a given $n$, its lower limit is $2(n-3)!$ when $c = 3$ and upper limit is $n!/3$ when $c = n$.

**Example 1**:

This problem is taken from [15] which was later adopted by Trandafir et al. [68]. It is a series system with 4 elements having the first and second order probabilities as:

$$[P_{ij}] = \begin{bmatrix} 0.27425312 & 0.17106964 & 0.13021655 & 0.09525911 \\ & 0.21185540 & 0.10920296 & 0.08120990 \\ & & 0.15865525 & 0.06566078 \\ & & & 0.11506967 \end{bmatrix}, \quad P_{ji} = P_{ij} \tag{21}$$

For notational convenience we have used $P_{ii} = P_i$ in Eq (4). Each element constitutes a minimal cut set in a series system and 4!=24 orderings of the minimal cut sets are possible for this problem. Figure 1 (left) shows the upper bound on $P[F_{sys}]$ for each of these orderings given by KVHD ($B_1$) and the proposed level 2 method ($B_2$). The relative errors (($B_1 - B_2$)/$B_1$) for all orderings are shown in Figure 1 (right). KVHD method yields its best $P_f = 0.363288$ for only 12 out of the 24 possibilities. Our level 2 method identifies every of those 12 cases, and an additional 6 orderings with the same best $P_f = 0.363288$. In each of the remaining six cases, our method improves upon KVHD. The lower scatter is evident from Figure 2: when all 24 orderings are considered, our level 2 upper bounds have a smaller mean (0.367) than KVHD bounds (0.379) and a significantly smaller coefficient of variation (COV=1.7%) than KVHD results (4.7%). Since the safety margins are jointly normal in the original problem statement, we can determine the exact system failure probability (0.349120) which is plotted as the horizontal line in Figure 1 (left).

While the level 2 bound in this example is clearly more effective than KVHD bound, we note that the best bound given by both are equal. We will come back to the question of whether the best bound can improve with increasing levels and if so under what conditions, but first, we present the general level *m* bound.

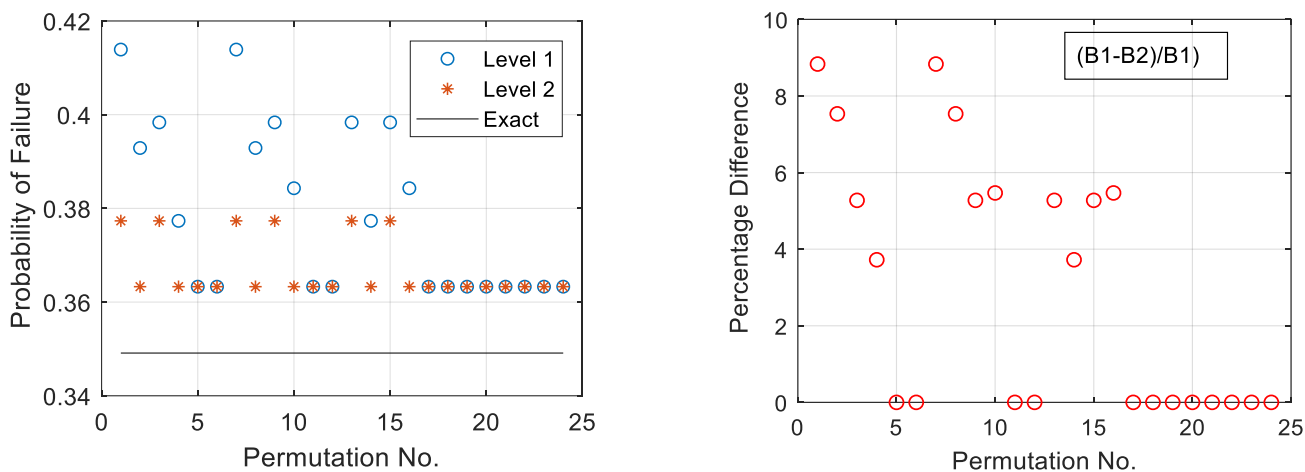

Figure 1: Four element series system: comparison of proposed level 2 with KVHD upper bound

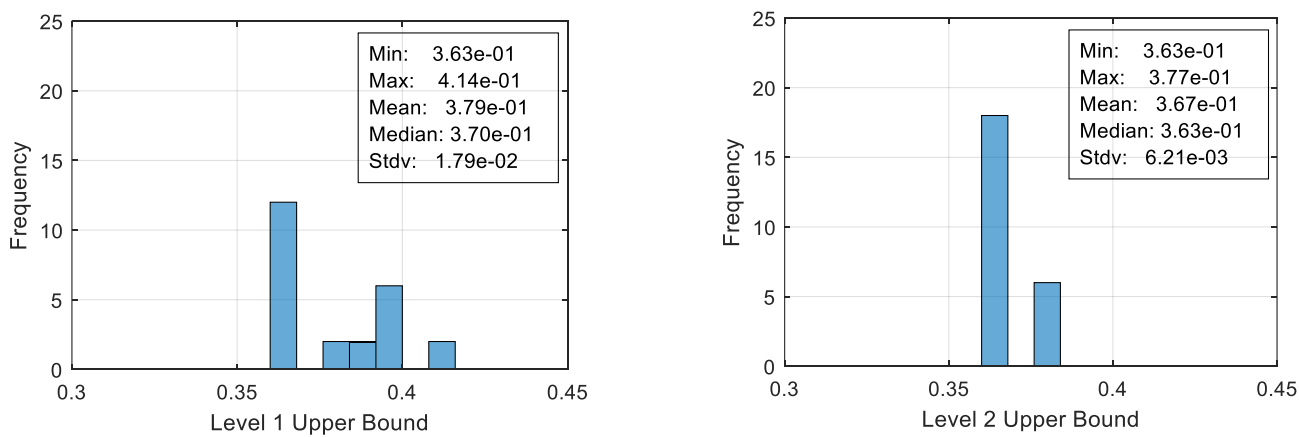

Figure 2: Four element series system: scatter in proposed level 2 vs. KVHD upper bound for all orderings of the index set

## 4. A nested hierarchy of upper bounds

The KVHD upper bound (7) and the upper bound derived in Eq (11) in fact belong to a hierarchy of second order bounds. KVHD bound considers only one second order intersection $C_{ij}$ in each line of Eq (6) whereas Eq (11) considers the union of two pairs $C_{ij}$ and $C_{jk}$ at a time. This bound can be further generalized by taking $m$ pairs at each line. To see this, take, for example, the union probability in the fourth line of Eq (6):

$$P^{(5)} = P(C_1C_5 \cup C_2C_5 \cup C_3C_5 \cup C_4C_5) \tag{22}$$

Since this term is subtracted, we need a lower bound to $P^{(5)}$ in order to derive an upper bound to $P[F_{sys}]$. For $m = 1$, that lower bound is simply the maximum of $\binom{4}{1} = 4$ terms, $\max_{j=1,\ldots,4}\{P_{j5}\}$:

$$P^{(5)} = P(C_1C_5 \cup C_2C_5 \cup C_3C_5 \cup C_4C_5) \geq \max_{j=1,\ldots,4}\{P_{j5}\} \tag{23}$$

For the level 2 bound, the lower bound to $P^{(5)}$ involves the maximum of $\binom{4}{2} = 6$ pair-wise union probabilities:

$$\begin{aligned}
P^{(5)} &= P(C_1C_5 \cup C_2C_5 \cup C_3C_5 \cup C_4C_5) \\
&\geq \max \begin{bmatrix} P(C_1C_5 \cup C_2C_5), P(C_1C_5 \cup C_3C_5), P(C_1C_5 \cup C_4C_5), \\ P(C_2C_5 \cup C_3C_5), P(C_2C_5 \cup C_4C_5), P(C_3C_5 \cup C_4C_5) \end{bmatrix} \\
&\geq \max_{\substack{1 \leq j,l < 5 \\ j \neq l}} \left[ P_{j5} + P_{l5} - \min(P_{j5}, P_{l5}, P_{lj}) \right]
\end{aligned} \tag{24}$$

Continuing this way, the lower bound to $P^{(5)}$ for $m = 3$ involves the maximum of $\binom{4}{3} = 4$ triplet-wise union probabilities as follows:

$$\begin{aligned}
P^{(5)} &= P(C_1C_5 \cup C_2C_5 \cup C_3C_5 \cup C_4C_5) \\
&\geq \max \begin{bmatrix} P(C_1C_5 \cup C_2C_5 \cup C_3C_5), P(C_1C_5 \cup C_2C_5 \cup C_4C_5) \\ P(C_1C_5 \cup C_3C_5 \cup C_4C_5), P(C_2C_5 \cup C_3C_5 \cup C_4C_5) \end{bmatrix} \\
&\geq \max_{\substack{1 \leq j,k,l < 5 \\ j \neq k,l \\ l \neq k}} \begin{bmatrix} P_{j5} + [P_{l5} - \min(P_{j5}, P_{l5}, P_{lj})] \\ +[P_{k5} - \min(P_{j5}, P_{k5}, P_{kj}) - \min(P_{l5}, P_{k5}, P_{kl})]^+ \end{bmatrix}
\end{aligned} \tag{25}$$

where $[a]^+ = \max[a,0]$. Generalizing, the level $m$ second order upper bound is:

$$P[F_{sys}] \leq \sum_{i=1}^{n}\left[ P_i - \max_{1\leq j_1 < j_2 < ... j_m < i}\left\{\sum_{r=1}^{m}\left[P_{j_r i} - \sum_{s=1}^{r-1}\min(P_{j_r i}, P_{j_s i}, P_{j_r j_s})\right]^+\right\}\right] = B_m, \quad m = 1,...,n-1 \quad (26)$$

which is the main result of this work. Eq (26) simplifies to Eq (7) for $m = 1$ and to Eq (11) for $m = 2$. By Theorem 2 we have shown that, given any permutation of the index set, the bound in Eq (26) for $m = 2$ is at least as good as that for $m = 1$. Here we generalize this to $m > 2$ as follows.

*Theorem* 3. In any ordering $(\pi)$ of the index set describing second order probabilities, the level $m + 1$ bound is less than or equal to the corresponding level $m$ bound, $m \leq n - 2$: $B_{m+1}(\pi) \leq B_m(\pi)$.

*Proof.* Incrementing $m$ by 1, we split the sum within the curly brackets of Eq (26) for any $1 \leq j_1, j_2, ..., j_m, j_{m+1} < i$ as,

$$\sum_{r=1}^{m+1}\left[P_{j_r i} - \sum_{s=1}^{r-1}\min(P_{j_r i}, P_{j_s i}, P_{j_r j_s})\right]^+ = \sum_{\substack{r=1 \\ r \neq v}}^{m+1}\left[P_{j_r i} - \sum_{s=1}^{r-1}\min(P_{j_r i}, P_{j_s i}, P_{j_r j_s})\right]^+ + \left[P_{j_v i} - \sum_{s=1}^{v-1}\min(P_{j_v i}, P_{j_s i}, P_{j_v j_s})\right]^+ \quad (27)$$

Since the second term on the RHS is non-negative,

$$\sum_{r=1}^{m+1}\left[P_{j_r i} - \sum_{s=1}^{r-1}\min(P_{j_r i}, P_{j_s i}, P_{j_r j_s})\right]^+ \geq \sum_{\substack{r=1 \\ r \neq v}}^{m+1}\left[P_{j_r i} - \sum_{s=1}^{r-1}\min(P_{j_r i}, P_{j_s i}, P_{j_r j_s})\right]^+ \quad (28)$$

Now taking maximum over all sequences $1 \leq j_1, j_2, ..., j_m, j_{m+1} < i$ and setting $v = m + 1$:

$$\max_{1\leq j_1, j_2,...j_m, j_{m+1} < i} \sum_{r=1}^{m+1}\left[P_{j_r i} - \sum_{s=1}^{r-1}\min(P_{j_r i}, P_{j_s i}, P_{j_r j_s})\right]^+ \geq \max_{1\leq j_1, j_2,...j_m < i} \sum_{r=1}^{m}\left[P_{j_r i} - \sum_{s=1}^{r-1}\min(P_{j_r i}, P_{j_s i}, P_{j_r j_s})\right]^+ \quad (29)$$

Subtracting both sides from $P_i$ and summing over $i = 1,...,n$ we get:

$$\sum_{i=1}^{n} P_i - \max_{1\leq j_1, j_2,...j_m, j_{m+1} < i} \sum_{r=1}^{m+1}\left[P_{j_r i} - \sum_{s=1}^{r-1}\min(P_{j_r i}, P_{j_s i}, P_{j_r j_s})\right]^+ \leq \sum_{i=1}^{n} P_i - \max_{1\leq j_1, j_2,...j_m < i} \sum_{r=1}^{m}\left[P_{j_r i} - \sum_{s=1}^{r-1}\min(P_{j_r i}, P_{j_s i}, P_{j_r j_s})\right]^+ \quad (30)$$

that is, the level $m + 1$ bound is at least as good as the level $m$ bound for any arbitrary permutation of the index set. Hence, proved.

We now generalize Condition 1 above and state Condition 2 under which the level $m + 1$ bound is strictly better than the level $m$ bound.

*Condition* 2: Given second order probabilities $P_{ji} = P_{ij}$, $i = 1,...,n-1$, $i < j \leq n$, in some ordering of the index set, $\{1, 2,...,n\}$, the terms satisfy

$$P_{j_r i} > \sum_{\substack{s=1 \\ s \neq r}}^{m+1} \min(P_{j_s i}, P_{j_r j_s}), \quad \forall r = 1, 2,..,m+1 < i, \text{ and } j_r, j_s < i \leq n \tag{31}$$

for every $\binom{i}{m+1}$ combination of the $m+1 < i$ indices.

It is easy to show that Condition 2 simplifies to Condition 1 for $m = 1$.

*Theorem* 4. If a particular ordering of the index set of second order probabilities satisfies Condition 2, the level $m + 1$ bound is less than the level $m$ bound in that ordering.

*Proof.* We have for one set of $m+1$ indices $1 \leq j_1, j_2,..., j_m, j_{m+1} < i \leq n$

$$\sum_{r=1}^{m+1}\left[P_{j_r i} - \sum_{s=1}^{r-1}\min(P_{j_r i}, P_{j_s i}, P_{j_r j_s})\right]^+ = \sum_{\substack{r=1 \\ r \neq v}}^{m+1}\left[P_{j_r i} - \sum_{s=1}^{r-1}\min(P_{j_r i}, P_{j_s i}, P_{j_r j_s})\right]^+ + \left[P_{j_v i} - \sum_{s=1}^{v-1}\min(P_{j_v i}, P_{j_s i}, P_{j_v j_s})\right]^+ \tag{32}$$

Since $P_{j_r i} > \sum_{\substack{s=1 \\ s \neq r}}^{m+1}\min(P_{j_s i}, P_{j_r j_s}) \Leftrightarrow P_{j_r i} > \sum_{\substack{s=1 \\ s \neq r}}^{m+1}\min(P_{j_r i}, P_{j_s i}, P_{j_r j_s})$ for each $r = 1, 2,..,m+1$, we have

$$\sum_{r=1}^{m+1}\left[P_{j_r i} - \sum_{s=1}^{r-1}\min(P_{j_r i}, P_{j_s i}, P_{j_r j_s})\right]^+ > \max_v(\sum_{\substack{r=1 \\ r \neq v}}^{m+1}\left[P_{j_r i} - \sum_{s=1}^{r-1}\min(P_{j_r i}, P_{j_s i}, P_{j_r j_s})\right]^+) \tag{33}$$

Now since this is true for every $m+1$ indices $1 \leq j_1, j_2,..., j_m, j_{m+1} < i \leq n$, we have

$$\max_{1 \leq j_1, j_2,..., j_m, j_{m+1} < i}\sum_{r=1}^{m+1}\left[P_{j_r i} - \sum_{s=1}^{r-1}\min(P_{j_r i}, P_{j_s i}, P_{j_r j_s})\right]^+ > \max_{1 \leq j_1, j_2,..., j_m < i}\sum_{r=1}^{m}\left[P_{j_r i} - \sum_{s=1}^{r-1}\min(P_{j_r i}, P_{j_s i}, P_{j_r j_s})\right]^+ \tag{34}$$

Subtracting both sides from $P_i$ and summing over $i = 1,...,n$ we get:

$$\sum_{i=1}^{n} P_i - \max_{1 \leq j_1, j_2, \ldots j_m, j_{m+1} < i} \sum_{r=1}^{m+1} \left[ P_{j_r i} - \sum_{s=1}^{r-1} \min(P_{j_r i}, P_{j_s i}, P_{j_r j_s}) \right]^+ <$$

$$\sum_{i=1}^{n} P_i - \max_{1 \leq j_1, j_2, \ldots j_m < i} \sum_{r=1}^{m} \left[ P_{j_r i} - \sum_{s=1}^{r-1} \min(P_{j_r i}, P_{j_s i}, P_{j_r j_s}) \right]^+ \quad (35)$$

i.e., $B_{m+1} < B_m$

Hence, proved.

**Example 1 (contd.):**

In Example 1 above, we find that Condition 2 is not satisfied in any ordering at level 2. Setting $i = 4$ and $m = 2$ in (31) and selecting $j_1 = 1, j_2 = 2, j_3 = 3$, Condition 2 requires,

$$\left. \begin{array}{l} P_{j_1 i} > \min(P_{j_2 i}, P_{j_1 j_2}) + \min(P_{j_3 i}, P_{j_1 j_3}) \\ P_{j_2 i} > \min(P_{j_1 i}, P_{j_1 j_2}) + \min(P_{j_3 i}, P_{j_2 j_3}) \\ P_{j_3 i} > \min(P_{j_1 i}, P_{j_1 j_3}) + \min(P_{j_2 i}, P_{j_2 j_3}) \end{array} \right\} \Rightarrow \left. \begin{array}{l} P_{14} > \min(P_{24}, P_{12}) + \min(P_{34}, P_{13}) \\ P_{24} > \min(P_{14}, P_{12}) + \min(P_{34}, P_{23}) \\ P_{34} > \min(P_{14}, P_{13}) + \min(P_{24}, P_{23}) \end{array} \right\} \quad (36)$$

Substituting the numerical values, we find the left hand sides of the three inequalities are respectively 0.09525911, 0.08120990, 0.06566078 while the right hand sides are:

$$\begin{array}{l} \min(0.08120990, 0.17106964) + \min(0.06566078, 0.13021655) = 0.14687068 \\ \min(0.09525911, 0.17106964) + \min(0.06566078, 0.10920296) = 0.16091989 \\ \min(0.09525911, 0.13021655) + \min(0.08120990, 0.10920296) = 0.22547566 \end{array} \quad (37)$$

It is straightforward to show that Condition 2 is not satisfied in every other permutation of the indices $j_1, j_2$ and $j_3$ as well. We can show the same to hold in every other ordering of the index set $\{1, 2, \ldots, n\}$ in this example.

**Example 2:**

In this example, we study how dependence among the cut sets affects the upper bounds. With $P_{ii} = P_i$, let the second order probabilities be of the form,

$$P_{ij} = P_i P_j + \delta, \quad i \neq j \quad (38)$$

The constant $\delta = 0$ if the cut sets $C_i$ and $C_j$ are statistically pairwise independent; if $\delta < 0$ the cut sets are negatively correlated and if $\delta > 0$ the cut sets are positively correlated. The allowable range of $\delta$ is:

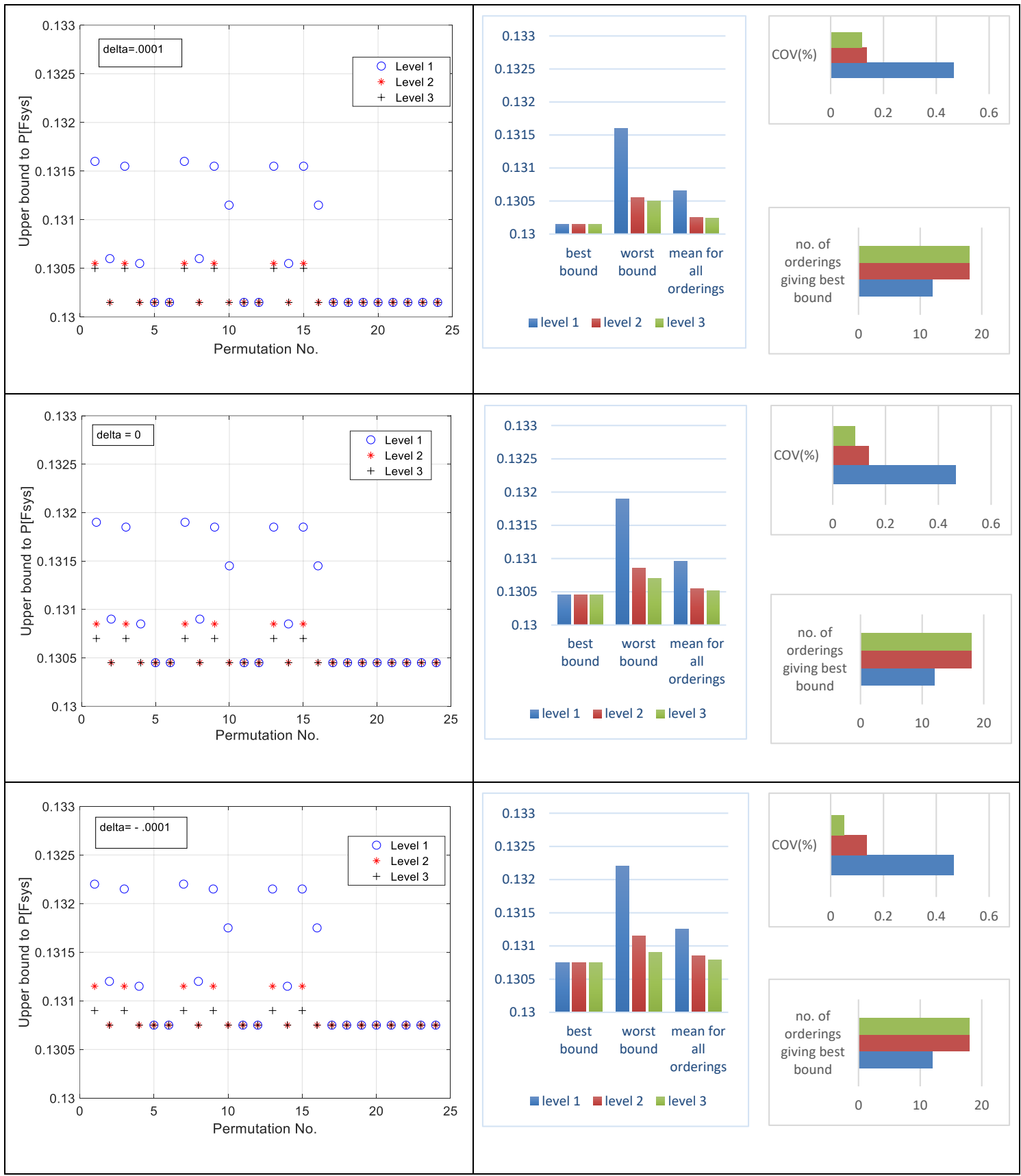

**Figure 3: Four element series system: comparison of levels 1 – 3 upper bounds: (a) top row - positively correlated element failures, (b) middle row – pairwise independent element failures, (c) bottom row – negatively correlated element failures**

$$-P_i P_j \leq \delta \leq \min(P_i, P_j) - P_i P_j \quad \forall i, j, i \neq j \tag{39}$$

We continue with a four element series system ($n = 4$), with first order failure probabilities $\{P_i\}$ = $[.01\ .025\ .03\ .07]^T$ and choose three value of $\delta \in \{0.0001, 0, -0.0001\}$, corresponding to positively correlated, pairwise independent and negatively correlated element failure events, respectively.

Figure 3 shows the levels 1, 2 and 3 bounds in all 24 permutations for each $\delta$. With $i = 4$, it is easy to check that Condition 2 is satisfied for $m = 1$ and $m = 2$ for all three values of $\delta$ in at least one ordering (i.e., $\{1, 2, 3, 4\}$) of the index set. In contrast to Example 1, we observe here the level 3 bound to be strictly better than the level 2 bound in 6 (and the level 1 bound in 12) out of 24 permutations of the index set, for each of the three cases of $\delta$. Thus, although the best (i.e., lowest) level 1, level 2 and level 3 upper bounds are all equal, level 1 achieves its best less frequently than do the higher levels. Further, the worst level 1 bound is significantly poorer than the worst level 2 bound, which in turn is significantly poorer than the worst level 3 bound. Further, when all 24 orderings are considered, the level 3 bounds show about 1/3 the scatter shown by level 2 bounds, and level 2 bounds in turn show about 1/3 the scatter shown by level 1 bounds.

**Example 3:**

We take a 5 element problem from [16]. The $5 \times 5$ second order symmetric probability matrix is:

$$[P_{ij}] = \begin{bmatrix} 4.548 & 1.776 & 1.790 & 1.559 & 0.119 \\ & 2.360 & 1.358 & 1.133 & 0.212 \\ & & 3.031 & 1.786 & 0.123 \\ & & & 2.744 & 0.269 \\ & & & & 1.469 \end{bmatrix} \times 0.01, \ P_{ji} = P_{ij} \tag{40}$$

$5! = 120$ permutations are possible for the index set and second order upper bounds up to the 4[th] level can be computed for each of those permutations. Table 1 lists a summary of the bounds. Clearly, levels 2 – 4 bounds are indistinguishable from one another, but level 1 bound performs significantly poorer than the higher level bounds: the level 1 bounds exhibit a much higher scatter, and the best level 1 bound equals the worst level 2 bound. Unlike level 1, the difference between the best and worst bounds at levels 2, 3 or 4 are insignificant. Because this is a small sized problem, the time taken to search through the 120 permutations are of the same order.

**Table 1: Summary of 120 upper bounds at 4 levels in Example 3**

|  | **Level 1** | **Level 2** | **Level 3** | **Level 4** |
|---|---|---|---|---|
| **Total CPU time (sec)** | 0.0132 | 0.0184 | 0.0647 | 0.0201 |
| **Minimum upper bound** | 0.08531 | 0.08438 | 0.08438 | 0.08438 |
| **Maximum upper bound** | 0.09241 | 0.08531 | 0.08531 | 0.08531 |
| **Mean upper bound** | 0.08847 | 0.08476 | 0.08476 | 0.08476 |
| **Median upper bound** | 0.08787 | 0.08442 | 0.08442 | 0.08442 |
| **COV (=SD/Mean) of upper bound (per cent)** | 2.52 | 0.53 | 0.53 | 0.53 |
| **Number of orderings giving minimum upper bound** | 12 | 12 | 12 | 12 |

**Example 4:**

This problem is taken from [9] as modified by Song and der Kiureghian [44]. A seven member determinate truss can fail due to the yielding of any of its seven members. Compression members are prevented from failing by buckling. The safety margins are:

$$M_i = X_i - L, \quad i=1,\dots,7 \tag{41}$$

The member strengths, $X_i$, are jointly normal: $X_1$ and $X_2$ each has a mean of 100kN and a standard deviation of 20kN while $X_3, \dots, X_7$ each has a mean of 200 kN and a standard deviation of 40 kN. The dependence structure is given by Dunnet–Sobel class correlation $\rho_{ij} = r_i r_j \ (i \neq j)$: $r_1 = 0.90$, $r_2 = 0.96$, $r_3 = 0.91$, $r_4 = 0.95$, $r_5 = 0.92$, $r_6 = 0.94$ and $r_7 = 0.93$ and $\rho_{ii} = 1$. The load $L = 100$kN is deterministic. The first order probabilities are all equal: $P_i = 1.88 \times 10^{-4}$. The complete second order probability matrix is:

$$[P_{ij}] = \begin{bmatrix} 18.8 & 5.73 & 4.35 & 5.42 & 4.59 & 5.13 & 4.85 \\ & 18.8 & 6.08 & 7.79 & 6.47 & 7.42 & 6.87 \\ & & 18.8 & 5.75 & 4.86 & 5.43 & 5.14 \\ & & & 18.8 & 6.10 & 6.88 & 6.48 \\ & & & & 18.8 & 5.76 & 5.44 \\ & & & & & 18.8 & 6.11 \\ & & & & & & 18.8 \end{bmatrix} \times 10^{-5}, P_{ji} = P_{ij} \qquad (42)$$

**Table 2: Summary of 5040 upper bounds at 6 levels in Example 4**

|  | Level 1 | Level 2 | Level 3 | Level 4 | Level 5 | Level 6 |
|---|---|---|---|---|---|---|
| **Total CPU time (sec)** | 0.0163 | 0.0333 | 5.586 | 7.739 | 9.601 | 9.049 |
| **Minimum upper bound** | 0.000912 | 0.000912 | 0.000912 | 0.000912 | 0.000912 | 0.000912 |
| **Maximum upper bound** | 0.000961 | 0.000944 | 0.000944 | 0.000944 | 0.000944 | 0.000944 |
| **Mean upper bound** | 0.000925 | 0.000919 | 0.000919 | 0.000919 | 0.000919 | 0.000919 |
| **Median upper bound** | 0.000924 | 0.000917 | 0.000917 | 0.000917 | 0.000917 | 0.000917 |
| **COV ( = SD/Mean) of upper bound (per cent)** | 1.22 | 0.83 | 0.83 | 0.83 | 0.83 | 0.83 |
| **Number of orderings giving minimum upper bound** | 24 | 1636 | 1636 | 1636 | 1636 | 1636 |

7!=5040 permutations of the minimal cut sets are possible for this problem. Multivariate normal integration yields the exact $P[F_{sys}]$ = 6.9988e−4. All levels give the lowest upper bound as 9.1216e−4: however the KVHD method yields this optimum for only 24 orderings, whereas the higher levels gives the lowest upper bound in almost a third of all cases (1636 out of 5040). Further, in 2420 non-optimal orderings, our method yields a smaller upper bound. The time taken, however, to search through the 5040 permutations is two orders of magnitude higher for levels 3 – 6 than for levels 1 and 2. As was the case with the two highest levels in Table 1, the time taken for the level 6 bound here is somewhat smaller than that for its preceding level because fewer terms need to be compared in the maximum value operation within the curly brackets of (26).

## 5. Does the optimal bound improve with levels?

The second order upper bound, for any level $m$, depends on the ordering of the index set. Let $B_m^*$ denote the best (i.e., smallest) level-$m$ bound $B_m$ identified across all orderings of the index set:

$$B_m^* = \min_{\substack{\text{all orderings } \pi \\ \text{of the index set}}} \left[ B_m(\pi) \right], \quad m = 1, \ldots, n-1 \tag{43}$$

We have shown that for *any* ordering of the index set, we must have $B_m(\pi) \geq B_{m+1}(\pi)$, that is, the level $m+1$ bound will always be as good or better than the level $m$ bound. We have also shown under what condition the relation becomes a strict inequality for a given ordering: $B_m(\pi) > B_{m+1}(\pi)$. Thus, while the first statement ensures that the optimal (i.e., best) bound over all orderings, $B_m^*$ in Eq (43), cannot get worse with increasing $m$, the second statement does not guarantee an improvement in the best. Additional conditions are required for $B_m^* > B_{m+1}^*$ to hold.

Without any loss of generality, let the second order probabilities, $P_{ij}$ ($i \neq j$), be all unique so that we can rank them as:

$$P^{[1]} > P^{[2]} > \ldots > P^{[n(n+1)/2]} \tag{44}$$

If some or all of them are equal, we can simply identify them interchangeably and the number of unique permutations will reduce. The best possible KVHD (i.e., level 1) upper bound is achieved if, for some ordering of the index set, the $i$th largest second order probability sits above the diagonal in column $i+1$ for each $i$. We denote such arrangements with the set $\pi^*$:

$$P^{[i]} = \max(P_{1,i+1}, P_{2,i+1}, \ldots, P_{i,i+1}; \pi^*) \tag{45}$$

which yields,

$$\begin{aligned} B_{1,\pi^*}^* &= P_1 + P_2 - P_{12} + \sum_{i=3}^{n} \left[ P_i - \max_{1 \leq j < i} \{P_{ji}\} \right] \\ &= \sum_{i=1}^{n} P_i - \sum_{i=1}^{n-1} P^{[i]} \end{aligned} \tag{46}$$

where the superscript '*' indicates the best possible value and $\pi^*$ refers to all those arrangements that satisfy (45).

We now look at the conditions necessary for the best level-2 bound to be better than the best level-1 bound, i.e., for $B_1^* > B_2^*$ to hold. For $n = 4$, the level 2 bound is:

$$B_2 = P_1 + P_2 - P_{12} + \sum_{i=3}^{4}\left[P_i - \max_{1 \le j < l < i}\left\{P_{ji} + P_{li} - \min(P_{ji}, P_{li}, P_{lj})\right\}\right]$$
$$= P_1 + P_2 - P_{12}$$
$$+ P_3 - \{P_{13} + P_{23} - \min(P_{13}, P_{23}, P_{12})\} \quad (47)$$
$$+ P_4 - \max\begin{cases} P_{14} + P_{24} - \min(P_{14}, P_{24}, P_{12}), P_{14} + P_{34} - \min(P_{14}, P_{34}, P_{13}), \\ P_{24} + P_{34} - \min(P_{24}, P_{34}, P_{23}) \end{cases}$$

The task is to place six second order probabilities above the diagonal of the probability matrix. We first restrict ourselves to Eq (45) since it ensures the best possible value of Ditlevsen's upper bound. Without any loss of generality we place the maximum $P^{[1]}$ among these at (1,2), then place $P^{[2]}$ in the third column and $P^{[3]}$ in the fourth column. A total of 2x3x3!=36 unique arrangements are possible involving $P^{[2]}, \ldots, P^{[6]}$ (another 36 arrangements can be made by interchanging the third and fourth columns; however these are not unique as they arise from a simple switching of the index set). Of these 36 arrangements, 20 show no improvement: $B_{1,\pi^*}^* = B_{2,\pi^*}^*$, another 4 yield $B_{1,\pi^*}^* > B_{2,\pi^*}^*$ conditionally, and the remaining 12 yield $B_{1,\pi^*}^* > B_{2,\pi^*}^*$ unconditionally. The cases are described in the following.

Let the indices $\{i, j, k\}$ be permutations of the integers $\{4, 5, 6\}$. Let $P^{[i]}$ be the other member in the third column (besides $P^{[2]}$). Thus $P^{[j]}$ and $P^{[k]}$ are elements of the fourth column.

[a] $P^{[2]}$ and $P^{[3]}$ are in different columns and in the same row (12 cases). If $P^{[i]} < \min(P^{[j]}, P^{[k]})$, i.e., $i = 6$, and if $P^{[4]} + P^{[5]} > P^{[3]} + P^{[6]}$ then $B_{1,\pi^*}^* > B_{2,\pi^*}^* = \Sigma P_i - P^{[1]} - P^{[2]} - [P^{[4]} + P^{[5]} - P^{[6]}]$ (2 cases). Otherwise, $B_{1,\pi^*}^* = B_{2,\pi^*}^* = \Sigma P_i - P^{[1]} - P^{[2]} - P^{[3]}$.

[b] $P^{[2]}$ and $P^{[3]}$ are in different columns and in different rows (24 cases). Let $P^{[2]}$ and $P^{[k]}$ be in the same row. If $P^{[i]} < P^{[j]}$, then $B_{1,\pi^*}^* > B_{2,\pi^*}^* = \Sigma P_i - P^{[1]} - P^{[2]} - [P^{[3]} + P^{[j]} - P^{[i]}]$ (2 cases). Otherwise, $B_{1,\pi^*}^* = B_{2,\pi^*}^* = \Sigma P_i - P^{[1]} - P^{[2]} - P^{[3]}$.

| $P_{i_1}$ | $P^{[1]}$ | $P^{[j_2]}$ | $P^{[j_3]}$ |
|---|---|---|---|
| | $P_{i_2}$ | $P^{[k_4]}$ | $P^{[k_5]}$ |
| Symmetric | | $P_{i_3}$ | $P^{[k_6]}$ |
| | | | $P_{i_4}$ |

| $P_{i_1}$ | $P^{[1]}$ | $P^{[k_4]}$ | $P^{[k_5]}$ |
|---|---|---|---|
| | $P_{i_2}$ | $P^{[j_2]}$ | $P^{[j_3]}$ |
| Symmetric | | $P_{i_3}$ | $P^{[k_6]}$ |
| | | | $P_{i_4}$ |

**Figure 4: Possible arrangements of the six unique second order probabilities in case [a] for 4x4 symmetric probability matrices. The diagonal terms are the first order probabilities and they can be placed without any restriction:** $(i_1, i_2, i_3, i_4)$ **are permutations of** (1,2,3,4). **The largest second order probability** $P^{[1]}$ **is placed at (1,2) without any loss of generality. In this case [a], the next two largest probabilities are in different columns but in the same row:** $(j_2, j_3)$ **are permutations of** (2,3). **The remaining three second order probabilities are placed in the remaining slots:** $(k_4, k_5, k_6)$ **are permutations of** (4,5,6).

The arrangements for case [a] are graphically shown in Figure 4. The other four cases can be depicted similarly. As stated above, identical results are obtained from 36 additional cases created by switching the third and fourth columns. We now relax the restriction imposed by Eq (45) and look at the remaining $2!\times 3!+3!\times 3! = 48$ cases (denoted by $\bar{\pi}$) where $P^{[2]}$ and $P^{[3]}$ are in the same column. Without any loss of generality, $P^{[1]}$ is still at (1,2). In $\bar{\pi}$, 12 arrangements show no improvement: $B^*_{1,\bar{\pi}} = B^*_{2,\bar{\pi}}$, another 4 yield $B^*_{1,\bar{\pi}} > B^*_{2,\bar{\pi}}$ conditionally, and the remaining 32 yield $B^*_{1,\bar{\pi}} > B^*_{2,\bar{\pi}}$ unconditionally. The cases are described in the following.

[c] $P^{[2]}$ and $P^{[3]}$ are in the 3rd column (12 cases). Regardless of where $P^{[i]}$, $P^{[j]}$ and $P^{[k]}$ are placed, there is no improvement: $B^*_{1,\bar{\pi}} = B^*_{2,\bar{\pi}} = \Sigma P_i - P^{[1]} - P^{[2]} - \max\{P^{[i]}, P^{[j]}, P^{[k]}\}$. Example 1 above is belongs to this case.

[d] $P^{[2]}$ and $P^{[3]}$ are in the 4th column and one of them is in (3,4). Then $B^*_{1,\bar{\pi}} > B^*_{2,\bar{\pi}}$ unconditionally (24 cases).

[e] $P^{[2]}$ and $P^{[3]}$ are in the 4th column and neither of them is in (3,4). Of the remaining terms with $i,j,k \in \{4,5,6\}$, let $P^{[i]}$ be the element in (3,4). If $P^{[i]} < \min(P^{[j]}, P^{[k]})$, i.e., $i = 6$, then there is no improvement: $B^*_{1,\bar{\pi}} = B^*_{2,\bar{\pi}} = \Sigma P_i - P^{[1]} - P^{[2]} - \max\{P^{[j]}, P^{[k]}\}$. If $P^{[i]} > \max(P^{[j]}, P^{[k]})$, i.e., $i = 4$, then there is definite improvement: $B^*_{1,\bar{\pi}} > B^*_{2,\bar{\pi}} = \Sigma P_i - P^{[1]} - P^{[2]} - \ldots\ldots$. Otherwise ($i = 5$), we have definite improvement ($B^*_{1,\bar{\pi}} > B^*_{2,\bar{\pi}}$) if $P^{[6]}$ is in the same row as $P^{[2]}$ and no improvement ($B^*_{\{1,\bar{\pi}\}} = B^*_{\{2,\bar{\pi}\}}$) if $P^{[6]}$ is not in the same row as $P^{[2]}$.

Combining the 120 results from arrangements $\pi$ and $\bar{\pi}$ described above, we find that 52 show no improvement, 56 show certain improvement, and the remaining 12 show improvement if certain conditions are satisfied. If the five probabilities are completely random, (i.e., $P_{ij} \sim U[0, \min(P_i, P_j)]$), the probability of finding $B^*_1 > B^*_2$ is $(56 + 4 \times 1/2 + 4 \times 1/2 + 4 \times 2/3)/120 = 52.2\%$ when $n = 4$.

We now show that this probability finding $B^*_1 > B^*_2$, provided the off-diagonal terms are conditionally independent and uniformly distributed, increases monotonically with $n$ and asymptotically reaches one.

*Theorem* 5. Given an $n$-dimensional matrix of second order probabilities $P_{ij}$ with IID diagonal elements $P_i \sim U[0,1]$ and conditionally independent off-diagonal elements $P_{ij} \sim U[0, \min(P_i, P_j)]$, the best level 2 bound is asymptotically better than the best level 1 bound: $\lim_{n \to \infty} P(B^*_2 < B^*_1) = 1$.

*Proof.* The $i$th lines in level 1 and level 2 bounds are, respectively, $P_i - L^1_i$ and $P_i - L^2_i$ where

$$L^1_i = \max_{j<i}(P_{ji})$$
$$L^2_i = \max_{j,k<i}\left(P_{ji} + P_{ki} - \min(P_{ji}, P_{ki}, P_{jk})\right)$$

(48)

It may be noted that,

$$L_i^1 \leq L_i^2, \ i \geq 3 \tag{49}$$

is always true and $L_i^1 = L_i^2$ for $i = 1$ and 2. Let $\pi_1^*$ be an ordering for which level 1 bound is optimal. We have already proved (Theorem 1) that for any ordering, the level 2 bound cannot be greater than the level 1 bound. Hence,

$$B_2(\pi_1^*) \leq B_1^* \tag{50}$$

Due to (49), B$_2$ is equal to $B_1^*$ if each of the line pairs $L_i^2, L_i^1$ are equal:

$$\{B_2(\pi_1^*) = B_1^*\} \Leftrightarrow (L_3^2 = L_3^1) \cap (L_4^2 = L_4^1) \cap \ldots (L_n^2 = L_n^1) \tag{51}$$

The complementary event gives the strict inequality,

$$\{B_2(\pi_1^*) < B_1^*\} \Leftrightarrow \{(L_3^2 = L_3^1) \cap (L_4^2 = L_4^1) \cap \ldots (L_n^2 = L_n^1)\}^c \tag{52}$$

Let us now consider the event,

$$T_i(\pi_1^*) = \{\min(P_{i-1}, P_i) U_i^1 > P_{i-1} U_i^2\} \cap \{\min(P_{i-2}, P_i) U_i^3 > P_{i-2} U_i^4\} \tag{53}$$

where $U_i^j \sim U(0,1), \ j = 1, \ldots, 4$ are independent standard uniform random variables. The probability of this event can be derived using an appropriate partition:

$$P[T_i(\pi_1^*)] = P[T_i(\pi_1^*) \cap P_i > P_{i-1} \cap P_i > P_{i-2}] + P[T_i(\pi_1^*) \cap P_i < P_{i-1} \cap P_i < P_{i-2}] + \\ P[T_i(\pi_1^*) \cap P_{i-1} > P_i > P_{i-2}] + P[T_i(\pi_1^*) \cap P_{i-1} < P_i < P_{i-2}] \tag{54}$$

The first term can be expanded as:

$$P[T_i(\pi_1^*) \cap P_i > P_{i-1} \cap P_i > P_{i-2}] = P[U_i^1 > U_i^2 \cap U_i^3 > U_i^4 \cap P_{i-1} < P_i \cap P_{i-2} < P_i] \\
= P[U_i^1 > U_i^2] P[U_i^3 > U_i^4] P[P_{i-1} < P_i \cap P_{i-2} < P_i] \\
= \frac{1}{2} \times \frac{1}{2} \times \int_0^1 \int_0^{p_i} \int_0^{p_i} dp_{i-2} \ dp_{i-1} \ dp_i = \frac{1}{2} \times \frac{1}{2} \times \frac{1}{3} = \frac{1}{12} \tag{55}$$

where we have used the mutual independence of $P_{i-2}, P_{i-1}, P_i$ and $U_i^j, j = 1, \ldots 4$. Proceeding similarly, the other three terms are, respectively, 1/54, 1/36 and 1/36, yielding the sum

$$P[T_i(\pi_1^*)] = \frac{1}{12} + \frac{1}{54} + \frac{1}{36} + \frac{1}{36} = \frac{17}{108} \tag{56}$$

Now, for any arbitrary quantity $P_b$, $T_i$ can be shown to be a subset of:

$$T_i(\pi_1^*) \subseteq \{\min(P_{i-1}, P_i) U_i^1 > \min(P_{i-1}, P_b) U_i^2\} \cap \{\min(P_{i-2}, P_i) U_i^3 > \min(P_{i-2}, P_{i-1}) U_i^4\} \tag{57}$$

which, using the definition given in the statement of this theorem, can be rewritten as:

$$T_i(\pi_1^*) \subseteq \{P_{i-1,i} > P_{i-1,b}\} \cap \{P_{i-2,i} > P_{i-2,i-1}\} \qquad (58)$$

Defining $P_{bi} = \max P_{ji}$, $j < i-1$, which implies $b \leq i-2$, the right hand side of (58) leads to:

$$\{P_{i-1,i} > P_{i-1,b}\} \cap \{P_{i-2,i} > P_{i-2,i-1}\}$$
$$\Rightarrow \{P_{i-1,i} + P_{bi} - P_{i-1,b} > P_{bi}\} \cap \{P_{i-2,i} + P_{i-1,i} - P_{i-2,i-1} > P_{i-1,i}\} \text{ where } b \leq i-2 \qquad (59)$$
$$\Rightarrow \{P_{i-1,i} + P_{bi} - \min(P_{i-1,i}, P_{i-1,b}, P_{bi}) > P_{bi}\} \cap \{P_{i-2,i} + P_{i-1,i} - \min(P_{i-2,i}, P_{i-2,i-1}, P_{i-1,i}) > P_{i-1,i}\}$$

Combining the LHS from both events gives a lower bound of the more general quantity $\max_{j,k<i, j\neq k} \{P_{ji} + P_{ki} - \min(P_{ji}, P_{jk}, P_{ki})\}$ while the combined RHS gives $\max_{j<i} P_{ji}$. In other words,

$$T_i(\pi_1^*) \subseteq \{L_i^2 > L_i^1\} \qquad (60)$$

which by (49) implies,

$$T_i(\pi_1^*) \subseteq \{L_i^2 = L_i^1\}^c \qquad (61)$$

The intersection of the complementary events, $T_i(\pi_1^*)^c$, has a probability bounded by:

$$P\left[\bigcap_{i=3}^n T_i(\pi_1^*)^c\right] \geq P\left[\bigcap_{i=3}^n \{L_i^2 = L_i^1\}\right] \text{ since } \bigcap_{i=3}^n T_i(\pi_1^*)^c \supseteq \bigcap_{i=3}^n \{L_i^2 = L_i^1\} \qquad (62)$$

Hence the probability of $B_2(\pi_1^*) < B_1^*$ in (52) can be bounded by:

$$P[B_2(\pi_1^*) < B_1^*] = 1 - P\left[\bigcap_{i=3}^n (L_i^2 = L_i^1)\right] \geq 1 - P\left[\bigcap_{i=3}^n T_i(\pi_1^*)^c\right] \qquad (63)$$

Since the events such as $T_3(\pi_1^*)$, $T_6(\pi_1^*)$, $T_9(\pi_1^*)$,... that are positioned at least 3 apart are mutually independent as they do not share any common elements, a lower bound to (63) can be obtained:

$$P[B_2(\pi_1^*) < B_1^*] \geq 1 - P\left[\bigcap_{i=3}^n T_i(\pi_1^*)^c\right] \geq 1 - P\left[\bigcap_{i=3,6,9,\ldots} T_i(\pi_1^*)^c\right] = 1 - \prod_{i=3,6,9,\ldots}\left(1 - P[T_i(\pi_1^*)]\right) \qquad (64)$$

Using the numerical value from (56),

$$P(B_2(\pi_1^*) < B_1^*) \geq 1 - \left(1 - \frac{17}{108}\right)^{\lfloor n/3 \rfloor} \qquad (65)$$

which, in the limit as the system size becomes large, yields

$$\lim_{n \to \infty} P(B_2(\pi_1^*) < B_1^*) = 1 \qquad (66)$$

Since $B_2(\pi_1^*)$ can only be greater than or equal to the level 2 optimum $B_2^*$, we must have

$$\lim_{n \to \infty} P\left(B_2^* < B_1^*\right) = 1 \tag{67}$$

Hence proved.

It can be shown that this asymptotic property holds for any two consecutive levels $m$ and $m + 1$, $1 \leq m \leq n-3$, with increasingly slower convergence. It can also be shown that for any finite $n$, the last two levels always have the same optimal bound: $B_{n-2}^* = B_{n-1}^*$.

Figure 5 shows the improvement in upper bounds from levels 1 through 4 with increasing system size in randomly generated second order probability matrices. Our level 2 bound is almost certain to show an improvement over KVHD bound as long as the second order probabilities are conditionally independent. The system has to be commensurately larger for higher level bounds to start showing noticeable improvements.

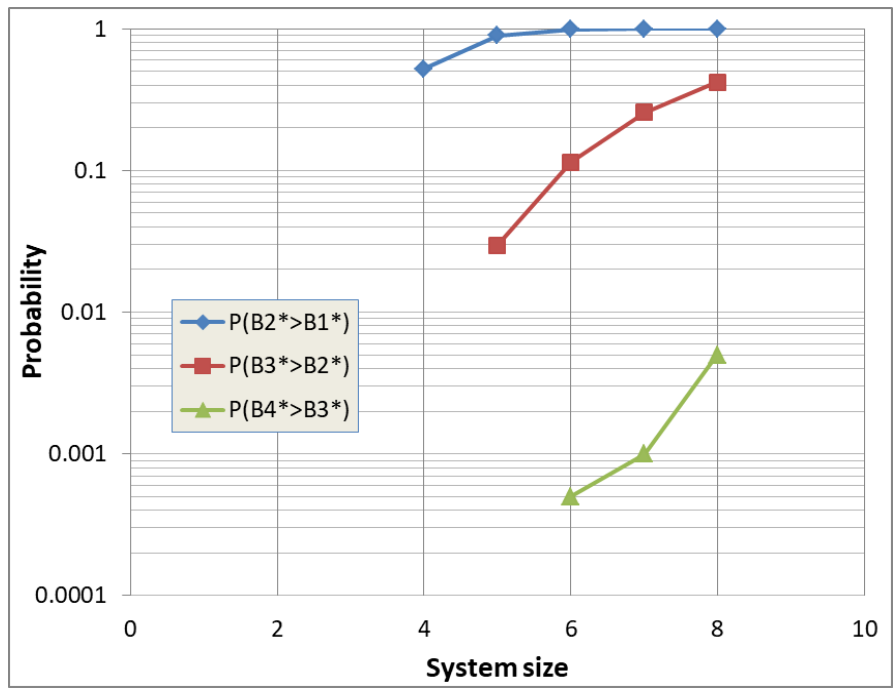

**Figure 5: Improvement in upper bounds from levels 1 through 4 with increasing system size in randomly generated second order probability matrices.**

**Example 5:**

In our final example, we look at one randomly generated 6x6 matrix used in Figure 5 which is reproduced as follows:

$$\left[P_{ij}\right] = \begin{bmatrix} 4.74467793 & 1.35693940 & 3.02042750 & 3.17568001 & 2.17177994 & 1.80796900 \\ & 2.34044502 & 0.58219757 & 0.38739530 & 0.19132633 & 1.39092307 \\ & & 3.60105675 & 0.44924975 & 0.33655831 & 1.88047290 \\ & & & 3.63910007 & 1.24586511 & 3.61723941 \\ & & & & 4.42818259 & 2.03204045 \\ & & & & & 6.94666654 \end{bmatrix} \times 10^{-3},$$

$$P_{ji} = P_{ij}$$

(68)

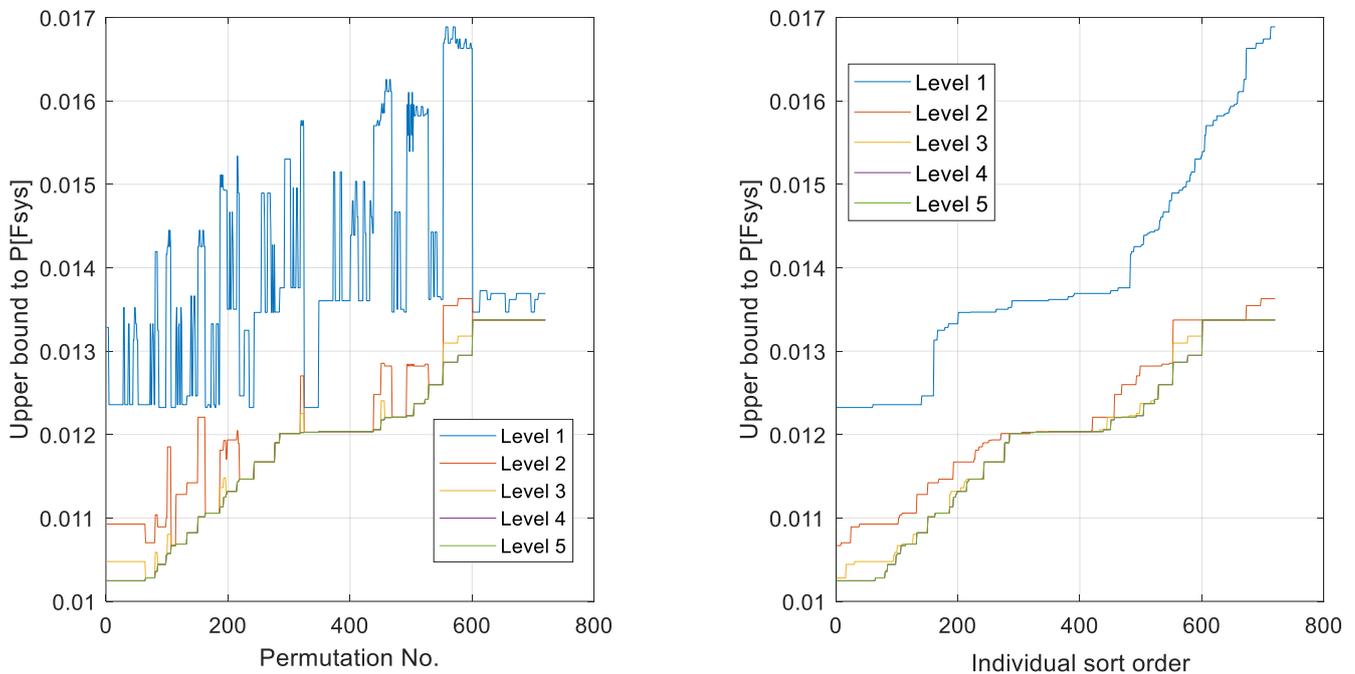

**Figure 6: Levels 1 – 5 upper bounds for one randomly generated 6x6 matrix used in Figure 5. Left: The 6! = 720 permutations of the index set are numbered by sorting the level 5 bound (green line) in increasing order. Since the bounds cannot worsen with increasing level, the five lines coincide segment-wise, but none of the 5 lines cross any other. The level 4 and level 5 bounds (green and purple lines) are coincident everywhere and, between the two, only the green line is visible. Right: each level is sorted individually and the values are presented in increasing order. Interestingly, the lines still do not cross each other. The starting point indicates the lowest possible value (i.e., B*) at each level. The best KVHD bound (.012324) is considerably larger than the best higher level bounds (0.010669, 0.010281, 0.010247 and .010247 respectively) although the benefit tapers off beyond level 3. At the other end, the worst value for each level presents a starker picture: KVHD bound performs much worse compared to the higher levels, and the higher level bounds stay confined within a noticeably narrow band.**

There are 6! permutations of the index set and Figure 6 (left) presents the five bounds corresponding to each of these 720 permutations: the permutations are numbered by sorting the level 5 bound (green line) in increasing order. By Theorem 3, the bounds cannot worsen with increasing level, and thus while they may coincide segment-wise, none of the 5 lines cross any other. The level 4 and level 5 bounds (green and purple lines) in this 6 × 6 problem are coincident everywhere and, between the two, only the green line is visible. The same results are presented differently in Figure 6 (right): each level is sorted individually and the values are presented in increasing order. It is interesting to note that the lines still do not cross each other. The starting point indicates the lowest possible value (i.e., $B_m^*$) at each level. The best KVHD bound (.012324) is considerably larger than the best higher level bounds (0.010669, 0.010281, 0.010247 and .010247 respectively) although the benefit tapers off beyond level 3. At the other end, the worst value for each level presents a starker picture: KVHD bound performs much worse compared to the higher levels, and the higher level bounds stay confined within a noticeably narrow band.

## 6. Conclusion

In this paper we derived a nested hierarchy of $m$-level second order upper bounds, $B_m$, on the union probability $P[F_{sys}] = P\left[\bigcup_{i=1}^{n} C_i\right]$ using only first and second order joint probabilities $P_i = P[C_i], P_{ij} = P[C_i C_j]$ since in practice, it is generally difficult to estimate joint probabilities beyond the second order. The well-known Kounias-Vanmarcke-Hunter-Ditlevsen (KVHD) bound - the current standard for upper bounds using second order joint probabilities – is the weakest member of this family ($m = 1$).

The tightness of such bounds depends on the particular ordering of the index set of the cut sets $C_i$ and identifying the optimal ordering is an important area of research. We proved that $B_m$ is non-increasing with level $m$ in every ordering of the cut sets, and derived conditions under which $B_{m+1}$ is strictly less than $B_m$ for any $m$ and any ordering. We also derived conditions under which the optimal (smallest, considering all $n!$ orderings of the index set) level $m + 1$ bound, $B_{m+1}^*$, is strictly less than the optimal level $m$ bound, $B_m^*$, and show that this improvement asymptotically achieves a probability of 1 as long as the second order joint probabilities are only constrained by the pair of corresponding first order probabilities but are otherwise independent.

Numerical examples showed that our second order upper bounds can yield tighter values than previously achieved, and in every case our bounds exhibit considerable less scatter across the *n*! permutations of the cut sets compared to KVHD bounds. Between successive levels, the highest relative improvement in the optimal $B_m^*$ for a given $n \times n$ second order probability matrix was found to occur between levels 1 and 2, and then to taper off at higher levels. The computation time increased with level *m*, however the increase from level 1 to level 2 is insignificant, which is also where the most improvement in $B_m^*$ is observed. Our results may lead to more efficient identification of the optimal upper bound when coupled with existing linear programming and tree search based approaches.

## Acknowledgment

The authors would like to thank the anonymous Reviewer 1 for pointing out the error in the originally submitted proof of Theorem 5 which has been rectified.

## References

1. Bonferroni CE. Teoria statistica classi e calcolo della proababilità. Pubbl R Ist Super Sci Econ Comm. 1936;8:1–63.
2. Chung KL, Erdos P. On the application of the Borel-Cantelli lemma. Transactions of the American Mathematical Society, American Mathematical Society. 1952;72:179-86.
3. Whittle P. Sur la distribution du maximum d'un polynome trigonomdtrique a coefficients aleatoires. Le Calcul des Probabilites et ses Applications, Colloques Internationaux de Centre National de la Recherche Scientifique. 1959;87:173-84.
4. Freudenthal AM, Garrelts JM, Shinozuka M. The analysis of structural safety. Journal of the Structural Division, ASCE. 1966;92(ST1):267–325.
5. Cornell CA. Bounds on the reliability of structural systems. Journal of the Structural Division, ASCE. 1967;93 (ST1):171–200.
6. Kounias EG. Bounds for the probability of a union with applications. The Annals of Mathematical Statistics, Institute of Mathematical Statistics. 1968;39: 2154-8.
7. Vanmarcke EH. Matrix formulation of reliability analysis and reliability-based design. Computers and Structures, Elsevier. 1973;3(4):757–70.
8. Hunter D. An upper bound for the probability of a union. Journal of Applied Probability, Cambridge University Press. 1976;13 (3):597–603.

9. Ditlevsen O. Narrow reliability bounds for structural systems. Journal of Structural Mechanics, Taylor and Francis. 1979;7(4):453–72.

10. Ahmed S, Koo B. Improved Reliability Bounds of Structural Systems. Journal of Structural Engineering, ASCE. 1990;116(11):3138-47.

11. Hohenbichler M, Rackwitz R. First-order concepts in systems reliability. Structural Safety, Elsevier. 1983;1(3):177–88.

12. Ramachandran K. Systems bounds: A critical study. Civil Engineering Systems, Taylor and Francis. 1984;1(3):123–8.

13. Feng Y. A method for computing structural system reliability with high accuracy. Computers & Structures, Elsevier. 1989;33(1):1–5.

14. Greig GL. An assessment of high-order bounds for structural reliability. Structural Safety, Elsevier. 1992;11(3-4):213–25.

15. Zhang YC. High-order reliability bounds for series systems and application to structural systems. Computers & Structures, Elsevier. 1993;46(2):381–6.

16. Ramachandran K. System reliability bounds: a new look with improvements. Civil Engineering and Environmental Systems, Taylor and Francis. 2004;21(4):265-78.

17. Cui W, Blockley D. On the bounds for structural system reliability. Structural Safety, Elsevier. 1991;9(4):247-59.

18. Qiu Z, Yang D, Elishakoff I. Probabilistic interval reliability of structural systems. International Journal of Solids and Structures, Elsevier. 2008;45(10):2850-60.

19. Wang C, Zhang H, Beer M. Computing tight bounds of structural reliability under imprecise probabilistic information. Computers & Structures, Elsevier. 2018;208:92-104.

20. Wang C, Zhang H, Beer M. Tightening the bound estimate of structural reliability under imprecise probability information. 13th International Conference on Applications of Statistics and Probability in Civil Engineering, ICASP13 May 26-30; Seoul, South Korea2019.

21. Zhang H, Mullen RL, Muhanna. Interval Monte Carlo methods for structural reliability. Structural Safety, Elsevier. 2010;32(3):183-90.

22. Hailerpin T. Best possible inequalities for the probability of a logical function of events. The American Mathematical Monthly, Mathematical Association of America. 1965;72(4):343–59.

23. Kounias S, Marin J. Best linear Bonferroni bounds. SIAM Journal on Applied Mathematics, Society for Industrial and Applied Mathematics. 1976;30(2):301–26.

24. Dawson OA, Sankoff D. An inequality for probabilities. Proceedings of the American Mathematical Society, American Mathematical Society. 1967;18 (3):504-7.

25. Gallot S. A bound for the maximum of a number of random variables. Journal of Applied Probability, Cambridge University Press. 1966;3(2):556-8.

26. Mallows CL. An inequality involving multinomial probabilities. Biometrika, Oxford University Press. 1968;55(2):422-4.

27. Sobel M, Uppuluri VRR. On Bonferroni-type inequalities of the same degree for the probability of unions and intersections. Annals of Mathematical Statistics, Institute of Mathematical Statistics. 1972; 43,(5):1549-58.

28. Kwerel SM. Most Stringent bounds on aggregated probabilities of partially specified dependent probability systems. Journal of the American Statistical Association, Taylor and Francis. 1975;70(350):472–9

29. Kwerel SM. Bounds on probability of a union and intersection of m events. Advances in Applied Probability, Cambridge University Press. 1975;7(2):431–48.

30. Galambos J. Bonferroni inequalities. Annals of Probability, Institute of Mathematical Statistics. 1977;5(4):577–81.

31. Galambos J, Mucci R. Inequalities for linear combinations of binomial moments. Publicationes Mathematicae Debrecen, University of Debrecen. 1980;27:263–8.

32. Platz O. A sharp upper probability bound for the occurrence of at least m out of n events. Journal of Applied Probability, Cambridge University Press. 1985;22(4):978–81.

33. Prekopa A. Boole-Bonferroni inequalities and linear programming. Operations Research, INFORMS. 1988;36(1):145–62

34. Prekopa A. Totally positive linear programming problems. Oxford Univ Press. 1989;L.V. Kantorovich Memorial Volume:197–207.

35. Prekopa A. Sharp bounds on probabilities using linear programming. Operations Research, INFORMS. 1990;38(2):227–39

36. Prekopa A. The discrete moment problem and linear programming. Discrete Applied Mathematics, Elsevier. 1990;27(3):235–54

37. Bukszár J, Prékopa A. Probability Bounds with Cherry Trees Mathematics of Operations Research, INFORMS. 2001;26(1):174-92.

38. Tomescu I. Hypertrees and Bonferroni inequalities. Journal of Combinatorical Theory, Elsevier. 1986;B 41:209–17.

39. Bukszár J, Szántai T. Probability bounds given by hypercherry trees. Optimization Methods and Software, Taylor and Francis. 2002;17(3):409-22.

40. Worsley KJ. An improved Bonferroni inequality and applications. Biometrika, Oxford University Press. 1982;69(2):297–302.

41. Boros E, Veneziani P. Bounds of degree 3 for the probability of the union of events. Piscataway, NJ Rutgers University, 2002.

42. Dohmen K. Bonferroni-type inequalities via chordal graphs. Combinatorics, Probability and Computing, Cambridge University Press. 2002;11:349–51.

43. Dohmen K. Lower bounds for the probability of a union via chordal graphs. arXiv: 10043416v2 [mathCO]. 2011.

44. Song J, Kiureghian AD. Bounds on system reliability by Linear Programming. Journal of Engineering Mechanics, ASCE. 2003;129(6):627–36.

45. Kiureghian AD, Song. J. Multi-scale reliability analysis and updating of complex systems by use of linear programming. Reliability Engineering and System Safety, Elsevier. 2008;93(2):288–97.

46. Chang Y, Mori Y. A Study of System Reliability Analysis Using Linear Programming. Journal of Asian Architecture and Building Engineering, Taylor and Francis. 2014;13(1):179-86.

47. Chang Y, Mai Y, Yi L, Yu L, Chen Y, Yang C, et al. Reliability Analysis of k-out-of-n Systems of Components with Potentially Brittle Behavior by Universal Generating Function and Linear Programming. Mathematical Problems in Engineering, Hindawi. 2020;2020(Article ID 8087242).

48. Byun JE, Song. J. Bounds on reliability of larger systems by linear programming with delayed column generation. Journal of Engineering Mechanics, ASCE. 2020;146(4):04020008-1 to -13.

49. J.Song, Kang WH, Lee YJ, Chun J. Structural System Reliability: Overview of Theories and Applications to Optimization. ASCE-ASME Journal of Risk and Uncertainty in Engineering Systems, Part A: Civil Engineering, ASCE. 2021;7(2): 03121001-1 to -24.

50. Prekopa A, Gao L. Bounding the probability of the union of events by aggregation and disaggregation in linear programs. Discrete Applied Mathematics, Elsevier. 2005;145(3):444–54.

51. de Caen D. A lower bound on the probability of the union. Discrete Mathematics, Elsevier. 1997;169(1-3):217–20.

52. Kuai H, Alajaji F, Takahara G. A lower bound on the probability of a finite union of events. Discrete Mathematics, Elsevier. 2000;215 147–58.

53. Yang J, Alajaji F, G.Takahara. On bounding the union probability using partial weighted information. Statistics & Probability Letters, Elsevier. 2016;116:38-44.

54. J.Yang, Alajaji F, Takahara G. Lower Bounds on the Probability of a Finite Union of Events. SIAM Journal on Discrete Mathematics, Society for Industrial and Applied Mathematics 2016;30(3):1437-52.

55. Feng C, Li L. On the Mori-Szekely conjectures for the Borel-Cantelli lemma. Studia Scientiarum Mathematicarum Hungarica, Akadémiai Kiadó. 2013;50:280–5.

56. Feng C, Li L, Shen J. Some inequalities in functional analysis, combinatorics, and probability theory. The Electronic Journal of Combinatorics. 2010;17:R58.

57. Mao Z, Cheng J, Shen J. A new lower bound on error probability for nonuniform signals over AWGN channels. Wireless Communications and Networking Conference (WCNC), IEEE; Piscataway, NJ2013. p. 3005–9.

58.	Cohen A, Merhav N. Lower bounds on the error probability of block codes based on improvements on de Caen's inequality. IEEE Transactions on Information Theory, IEEE Information Theory Society. 2004;50(2):290-310.

59.	Szántai T. Improved Bounds and Simulation Procedures on the Value of the Multivariate Normal Probability Distribution Function. Annals of Operations Research, Springer. 2000;100:85–101

60.	Georgakopoulos G, Kavvadias D, Papadimitriou CH. Probabilistic Satisfiability Journal of Complexity, Elsevier. 1988;4(1):1-11.

61.	Boros E, Scozzari A, Tardella F, Veneziani P. Polynomially computable bounds for the probability of the union of events. Mathematics of Operations Research, INFORMS. 2014;39(4):1311–29.

62.	Zemel E. Polynomial algorithms for estimating network reliability. Networks, Wiley. 1982;12(4):439–52.

63.	Jaumard B, Hansen P, de Aragão MP. Column generation methods for probabilistic logic. INFORMS Journal on Computing, INFORMS. 1991;3(2):135–48.

64.	MM MMD, Laurent M. Geometry of Cuts and Metrics Springer-Verlag, Berlin; 1997.

65.	Boros E, Hammer PL. Cut-polytopes, Boolean quadric polytopes and nonnegative quadratic pseudo-Boolean functions. Mathematics of Operations Research, INFORMS 1993;18(1):245–53.

66.	Kavvadias D, Papadimitriou CH. A linear programming approach to reasoning about probabilities. Annals of Mathematics and Artificial Intelligence, Springer. 1990;1:189–205.

67.	Veneziani P. Upper bounds of degree 3 for the probability of the union of events via linear programming. Discrete Applied Mathematics, Elsevier 2009;157(4):858–63.

68.	Trandafir R, Demetriu S, Mazilu IM. Determination of reliability bounds for structural systems using linear programming. Analele Ştiinţifice ale Universităţii "Ovidius. 2003;11(2):171-82.